\documentclass[12pt,a4paper]{article}

\usepackage{latexsym}
\usepackage{amsfonts}
\usepackage{amssymb}

\newcommand{\be}{\begin{equation}}
\newcommand{\ee}{\end{equation}}
\newcommand{\bea}{\begin{eqnarray}}
\newcommand{\eea}{\end{eqnarray}}

\def\G{{\cal G}}                                 %
\def\H{{\cal H}}                                 %
\def\F{{\cal F}}                                 %
\def\A{{\cal A}}                                 %
\def\M{{\cal M}}                                 %
\def\K{{\cal K}}                                 %
\def\cR{{\cal R}}                                %
\def\cL{{\cal L}}                                %
\def\Ad{{\mathrm{Ad}}}                           
\begin{document}

\vspace*{0.5cm}
\begin{center}
{\Large \bf Poisson-Lie dynamical $r$-matrices  from Dirac reduction}
\end{center}

\vspace{1.0cm}

\begin{center}
L. Feh\'er\footnote{Postal address: MTA KFKI RMKI,
H-1525 Budapest 114, P.O.B. 49,  Hungary} \\

\bigskip

{\em Department of Theoretical Physics\\
MTA  KFKI RMKI and University of Szeged \\
E-mail: lfeher@rmki.kfki.hu\\ }

\end{center}

\vspace{1.0cm}

\begin{abstract}

The Dirac reduction technique
used previously to obtain solutions of the classical dynamical
Yang-Baxter equation on the dual of a Lie algebra is extended
to the Poisson-Lie case and is shown to yield naturally
certain dynamical $r$-matrices on the duals
of Poisson-Lie groups found by Etingof, Enriquez and Marshall
in math.QA/0403283.

\end{abstract}

\bigskip



\newpage

\section{Introduction}

Let $G$ be a connected Poisson-Lie (PL) group of
the coboundary type, and denote by $\G:=\mathrm{Lie}(G)$
its Lie algebra.
The Poisson bracket (PB) on $G$ can be encoded by the formula
\be
\{ g_1, g_2\}_G = [g_1 g_2,R], \qquad g\in G,
\label{1.1}\ee
where $R\in \G\wedge \G$ solves the (modified) classical
Yang-Baxter equation
\be
[R_{12}, R_{13}]+[R_{12}, R_{23}] +
[R_{13}, R_{23}] ={\cal I}_R,
\label{1.2}\ee
with some $G$-invariant ${\cal I}_R\in \G\wedge\G\wedge\G$.
Consider a PL subgroup  $K\subseteq G$ and the corresponding
dual PL  group $K^*$.
Fix some open submanifold $\check K^*\subset K^*$.
By definition, a {\em PL dynamical $r$-matrix with respect to the pair}
$K\subseteq G$ is an `admissible'  (smooth or meromorphic) mapping
$r: \check K^* \rightarrow \G\wedge \G$, which is $\K$-equivariant
in the natural sense and satisfies the  equation
\be
[R_{12} + r_{12}, R_{23} + r_{23}] + K_1^a \cL_{K_a} r_{23}
+ \mathrm{cycl. perm.} = {\cal I}_{R,r},
\label{1.3}\ee
where $\{ K^a\} \subset \K:=\mathrm{Lie}(K)$,
$\{K_a\} \subset \K^*=\mathrm{Lie}(K^*)$ are bases in duality,
${\cal I}_{R,r} \in \G\wedge \G\wedge \G$ is a $G$-invariant constant, and
$\cL_{K_a}$ is the left-derivative associated with $K_a\in \K^*$.
Equation (\ref{1.3}) is called the PL-CDYBE
 for the pair $K\subseteq G$.
(The shorthand CDYBE stands for `classical
dynamical Yang-Baxter equation'.)
Motivated by the study of the PL symmetries of the chiral
Wess-Zumino-Novikov-Witten  phase space \cite{BFP},
this equation was considered for $K=G$ in \cite{FM1,FM2}.
The general case of a proper PL subgroup $K\subset G$
(also without restricting $G$ to be of the coboundary type)
was investigated in \cite{EEM}.
See also \cite{DM, M} for an even more general notion of PL
dynamical $r$-matrices.
 If $R$ is set to zero, then any Lie subgroup $K\subset G$ is a PL subgroup,
and the dual $K^*$ becomes $\K^*$ with its
linear Lie-Poisson structure.
Thus for $R=0$ the PL-CDYBE reproduces the CDYBE for the pair
$\K\subseteq \G$ as defined in \cite{EV}.

Etingof and Varchenko \cite{EV} introduced
a useful technique of reduction of variables that
connects, for example, solutions of the CDYBE on
$\G^*$ and on $\K^*$, where $\K$ is a Levi subalgebra of a simple
Lie algebra $\G$.
In \cite{FGP} the reduction
technique of \cite{EV} was shown to be equivalent to
the application of a suitable
Dirac reduction to the PL groupoid that underlies
the geometric interpretation of the CDYBE.
The reduction technique of \cite{EV} (see also \cite{EE})
has been generalized in \cite{EEM} to the PL case,
leading to new PL dynamical $r$-matrices.

The purpose of the present note is to show that, as was anticipated
in \cite{FM1}, the Dirac
reduction method of \cite{FGP} extends naturally to the PL case, too.
This method permits us to obtain a better understanding
of some constructions in  \cite{EEM}, and
it may prove useful in future investigations as well.

For simplicity, we shall focus on the {\em triangular}
PL dynamical $r$-matrices, which satisfy the extra
condition ${\cal I}_{R,r}={\cal I}_R$ by definition.
In our Dirac reduction the starting phase space will be
the manifold
$G \times \check K^*$
equipped with a PB encoding a triangular
PL dynamical $r$-matrix $r: \check K^* \rightarrow \G\wedge \G$.
If $H\subset K$ is a PL subgroup and certain further conditions
are satisfied, then Dirac reduction yields
$G\times \check H^*$
in such a way that the reduced PB (the `Dirac bracket')
encodes another triangular PL dynamical $r$-matrix
$r^*: \check H^* \rightarrow \G\wedge \G$.
In particular, if the starting $r$-matrix is zero, then
we recover the $\sigma_\H^\G$ family of $r$-matrices discovered in
\cite{EEM}.
The conditions guaranteeing our Dirac reduction to work
are also used in the direct definition of $\sigma_\H^\G$
in \cite{EEM}.
Here the conditions will be seen to emerge
naturally from the construction.

The main part of the paper is Section 3, where we deal with
the Dirac reduction of $K^*$ to $\check H^*$ and its application
to the triangular solutions of (\ref{1.3}).
The general (non-triangular) case, together with
examples and their possible applications, is briefly
discussed in Section 4.

\section{Geometric model for triangular dynamical $r$-matrices}
 \setcounter{equation}{0}

The PL-CDYBE is encoded by the Jacobi identities of the PB
on certain Poisson manifolds.
The following model, valid
in the triangular case, can be found in \cite{EEM}.

Below, for any Lie group $A$ the adjoint action of $a\in A$
on $X\in \A:=\mathrm{Lie}(A)$ is denoted simply by $\Ad_a(X)=aX a^{-1}$.
In the same spirit, regarding $a\in A$ as a matrix,
we may write the left and right derivatives as  $\cL_X a= Xa$, $\cR_X a= aX$ and so on.
For an arbitrary function $f$ on $A$, we have
$({\cL}_X f)(a):= \frac{d}{d t} f(e^{t X} a)\Big\vert_{t=0}$ and
$({\cR}_X f)(a)$ is defined similarly.
Correspondingly, the $\A^*$-valued left and right `gradients' are
defined by
\be
 \langle \nabla_a f, X \rangle =(\cL_X f)(a),
\qquad
 \langle \nabla'_a f, X \rangle=(\cR_X f)(a).
\label{2.1}\ee
We need to recall (see, e.g., \cite{LW}) that the PB on the dual $K^*$ of a PL group $K$
can be written as
\be
\{ f_1, f_2 \}_{K^*}(\kappa)= \langle\langle \nabla_\kappa f_1,
\kappa (\nabla'_\kappa f_2) \kappa^{-1} \rangle\rangle.
\label{2.2}\ee
Here $\nabla_\kappa f_i, \nabla'_\kappa f_i\in \K=(\K^*)^*$,
$\langle\langle\ ,\ \rangle\rangle$ is the `scalar product'
on the Drinfeld double Lie algebra ${\cal D}({\K,\K^*})$, the
adjoint action of $\kappa \in K^*$ on $X\in \K$ refers
to the Drinfeld double Lie group $D({K,K^*})$ that contains $K$
and $K^*$ as Lie subgroups.
We need also the infinitesimal left dressing action
of $\K$ on $K^*$, which is defined by the formula
\be
\mathrm{dress}_X \kappa = \kappa (\kappa^{-1} X \kappa)_{\K^*},
\qquad \forall X\in \K,
\label{2.3}\ee
where we use the decomposition of $\forall Y\in {\cal D}({\K,\K^*})$ into
$Y=Y_\K + Y_{\K^*}$ with $Y_\K\in \K$, $Y_{\K^*}\in \K^*$.

Fixing $\check K^*$ to be an open submanifold
of $K^*$, consider the manifold
\be
Q:= G \times \check K^* = \{ (g, \kappa)\}.
\label{2.4}\ee
We write $Q(\check K^*)$ if we want to emphasize the
dependence on the choice of $\check K^*$.
For functions $\phi$ on $G$ and $f$ on $\check K^*$,
let $\phi'$ and $\hat f$ be the functions on $Q$ given by
\be
\phi'(g,\kappa) =\phi(g),
\qquad
\hat f(g, \kappa)=f( \kappa).
\label{2.5}\ee
Take an admissible function
\be
r:\check K^*\rightarrow \G\wedge\G
\label{2.6}\ee
and try to define a PB on $Q$ by means of the ansatz
\bea
&\{ \hat f_1, \hat f_2\}_{Q}( g,  \kappa)=
\{ f_1, f_2 \}_{\check K^*}(\kappa),
\label{2.7}\\
&\{ \phi', \hat f\}_{Q}( g, \kappa)=
\langle \nabla'_g \phi, \nabla_{\kappa} f\rangle,
\label{2.8}\\
&\{ \phi'_1, \phi'_2 \}_{ Q}( g,\kappa)=
\langle \nabla'_g \phi_1 \otimes \nabla'_g \phi_2, R+ r( \kappa) \rangle
- \langle \nabla_g \phi_1 \otimes \nabla_g \phi_2, R \rangle,
\label{2.9}\eea
where $R$ is the underlying constant solution of (\ref{1.2}).
One can verify

\medskip
\noindent
{\bf Proposition 2.1.}
{\em The bracket $\{\ ,\ \}_{Q}$ satisfies
the Jacobi identity if and only if the infinitesimal
equivariance condition
\be
\mathrm{dress}_X r = [X\otimes 1 + 1\otimes X, r]
\qquad \forall X\in \K,
\label{2.10}\ee
and the PL-CDYBE (\ref{1.3})
with ${\cal I}_{R,r}={\cal I}_R$ hold for $r$.}

\section{Dirac reduction and dynamical $r$-matrices}
 \setcounter{equation}{0}

Let $H\subset K \subseteq G$ be a chain of connected PL subgroups of $G$.
Given a Poisson manifold
$(Q(\check K^*), \{\ ,\ \}_{Q(\check K^*)})$,
we  wish to reduce it to a Poisson manifold of the same kind,
but with respect to the subgroup $H\subset K$.
We wish to achieve this by viewing $Q(\check H^*)$ as a submanifold
of $Q(\check K^*)$ specified by {\em second class}
constraints in Dirac's sense \cite{Dirac}. Crucially, the constraints must be
such that the reduced PB (the `Dirac bracket') resulting from
$\{\ ,\ \}_{Q(\check K^*)}$ should have the form of
$\{\ ,\ \}_{Q(\check H^*)}$. If this happens, then
the triangular $r$-matrix $r: \check K^* \rightarrow \G\wedge \G$
contained in $\{\ ,\ \}_{Q(\check K^*)}$
gives rise to a reduced triangular $r$-matrix
$r^*: \check H^* \rightarrow \G\wedge \G$ contained in
$\{\ ,\ \}_{Q(\check H^*)}$.

\subsection{Dirac reduction of $K^*$ to $\check H^*$}

The reduction of $r$-matrices sketched above can only work
if an open submanifold of $( H^*, \{\ ,\ \}_{H^*})$
can be obtained as the Dirac reduction of
$(K^*, \{\ ,\ \}_{K^*})$.
To investigate the condition for this, let
${\cal D}(\K,\K^*)$ and ${\cal D}(\H,\H^*)$ be the Drinfeld doubles of the
Lie bialgebras corresponding to the PL groups $K$ and $H$.
As linear spaces,
\be
{\cal D}(\K,\K^*)=\K + \K^*,
\qquad
{\cal D}(\H,\H^*)=\H+\H^*,
\label{3.1}\ee
where
$\H=\mathrm{Lie}(H)$, $\H^*=\mathrm{Lie}(H^*)$ and
similarly for $K$.

We have assumed that $H\subset K$ is a connected PL subgroup,
and this is known \cite{STS} to be equivalent to the condition that
$\H^\perp\subset \K^*$,
\be
\H^\perp = \{ \alpha\in \K^*\,\vert\, \langle \alpha, X\rangle =0\quad
\forall X\in \H\,\},
\label{3.2}\ee
is an ideal of the Lie subalgebra $\K^*\subset {\cal D}(\K,\K^*)$.
Next, $H^*$ must clearly be a Lie subgroup of $K^*$ for our
construction, and this requires that
$\H^* \subset \K^*$ must be a Lie subalgebra.
We can encode these data in a vector space decomposition
\be
\K=\H + \M,
\label{3.3}\ee
which induces
\be
\K^*=\H^* + \M^*
\qquad\hbox{with}\qquad
\H^*=\M^\perp, \quad
\M^*=\H^\perp.
\label{3.4}\ee
In addition to $\H^*$ being a Lie subalgebra and $\M^*$ being a Lie
ideal, we shall need (see also Remark 3.5 below) the decomposition (\ref{3.3}) to be reductive
\be
[\H,\M]\subset \M,
\label{3.5}\ee
and of course the constraints specifying $\check H^*$ inside
$K^*$ must be second class.
Let $\{ M^i\}\subset \M$ be a basis.
The second class nature of the constraints turns out equivalent to
the non-degeneracy of the matrix
\be
C^{ij}(\lambda)= \langle\langle (\lambda M^i \lambda^{-1})_{\M},
\lambda M^j \lambda^{-1}\rangle\rangle
\quad \hbox{for}\quad \lambda \in \check H^*,
\label{3.6}\ee
defined using the Drinfeld double $D(K,K^*)$.

We next show that (\ref{3.6})  together with the foregoing other
assumptions guarantees
the desired reduction of $K^*$ to $\check H^*$.
We begin by proving some auxiliary statements.

\medskip
\noindent
{\bf Lemma 3.1.}
{\em With the above notations, suppose that
$\H \subset \K$ and $\H^*\subset \K^*$ are
Lie subalgebras, $\H^\perp \subset \K^*$
is a Lie ideal and $[\H,\M]\subset \M$.
Then $\H+\H^*$ is a Lie subalgebra of the double ${\cal D}(\K,\K^*)$.
This subalgebra of ${\cal D}(\K,\K^*)$ can be identified with
the double ${\cal D}(\H,\H^*)$.}

\medskip
\noindent
{\em Proof.}
We need to show that $[\H,\H^*]\subset \H+\H^*$ inside ${\cal D}(\K,\K^*)$.
With the Lie bracket $[\ ,\ ]$ and invariant `scalar product'
$\langle\langle\ ,\ \rangle\rangle$ of ${\cal D}(\K,\K^*)$, we have
$$
\langle\langle  [\H^*,\H],\M\rangle\rangle =
\langle\langle  \H^*, [\H,\M]\rangle\rangle
\subset
\langle\langle  \H^*, \M\rangle\rangle = \{ 0\},
$$
since $[\H,\M]\subset \M$, and
$$
\langle\langle  [\H,\H^*],\M^* \rangle\rangle =
\langle\langle  \H, [\H^*,\M^*]\rangle\rangle
\subset
\langle\langle  \H, \M^*\rangle\rangle = \{ 0\},
$$
since $[\H^*,\M^*]\subset \M^*$ as $\M^*=\H^\perp$.
{\em Q.E.D.}

\medskip
\noindent
{\bf Lemma 3.2.}
{\em
Under the assumptions of Lemma 3.1,
consider the connected Lie subgroup $H^* \subset K^*$ corresponding to
$\H^* \subset \K^*$.
Parametrize the elements in some neighbourhood of $H^*$ in $K^*$
as
\be
\kappa = \lambda e^\mu \qquad \lambda\in H^*,\quad \mu\in \check \M^*,
\label{3.9}\ee
where $\check \M^*$ is some neighbourhood of zero in $\M^*$.
Take a function $F\in \F(H^*)$ and extend it (locally) to
$f\in \F(K^*)$ by
\be
f(\lambda e^\mu) = F(\lambda).
\label{3.10}\ee
Then
\be
\nabla_\lambda f=\nabla_\lambda F
\quad\hbox{and}\quad
\nabla'_\lambda f =\nabla_\lambda' F.
\label{3.11}\ee
}

\medskip
\noindent
{\em Proof.}  In principle, $\nabla_\lambda f\in \K$ and
$\nabla_\lambda F\in \H$.
For $X\in \H^*$ we have
$f(e^{tX} \lambda)=F(e^{tX} \lambda)$, and for $Y\in \M^*$ we have
$f(e^{tY}\lambda)= f(\lambda \lambda^{-1} e^{tY} \lambda)
= f(\lambda e^{t \lambda^{-1} Y \lambda}) = F(\lambda)$ since
$\lambda^{-1} Y \lambda \in \M^*$ by $[\H^*,\M^*]\subset \M^*$.
This implies the first equality in (\ref{3.11}).
The second equality follows
similarly, and actually it is also a consequence of the first one.
Indeed, $\nabla'_\lambda F= (\lambda^{-1} \nabla_\lambda F\lambda)_\H$
in the double of $\H$, and
$\nabla'_\lambda f= (\lambda^{-1} \nabla_\lambda f\lambda)_\K$
in the double of $\K$ on general grounds, which implies the second
equality by Lemma 3.1.
{\em Q.E.D.}

\medskip
\noindent
{\bf Lemma 3.3.}
{\em
Keeping the preceding assumptions,
for a constant $M\in \M$ define the function $\xi_M$ on
a neighbourhood of $H^*$ in $K^*$ by
\be
\xi_M(\lambda e^\mu):= \langle \mu, M \rangle
\qquad (\lambda\in H^*,\,\mu\in \check \M^*).
\label{3.12}\ee
For this function, we have
\be
\nabla'_\lambda \xi_M = M,
\qquad
\nabla_\lambda \xi_M = (\lambda M \lambda^{-1})_\K=
(\lambda M \lambda^{-1})_\M.
\label{3.13}\ee
As a consequence,
\be
\{ f, \xi_M\}_{K^*}(\lambda)=0 \qquad (\forall \lambda\in H^*)
\label{3.14}\ee
for any functions $f$ and $\xi_M$ defined in (\ref{3.10}), (\ref{3.12}).}

\medskip
\noindent
{\em Proof.}
It is simple to confirm
$\nabla'\xi_M =M$ directly from the definition,
and
this implies
$\nabla_\lambda \xi_M= (\lambda M \lambda^{-1})_\K$ by the
universal connection between
left and right derivatives.
The last equality in (\ref{3.13}) follows since
$\langle\langle \lambda \M \lambda^{-1}, \H^*\rangle\rangle =
\langle\langle  \M, \lambda^{-1} \H^* \lambda \rangle\rangle
\subset \langle\langle \M, \H^*\rangle\rangle = \{0\}$.
By using (\ref{2.2}),
the statement of (\ref{3.14}) is a consequence of the fact that
$\nabla_\lambda f\in \H$ and $\nabla'_\lambda \xi_M \in \M$.
Indeed,
$\langle\langle \H, \lambda \M \lambda^{-1} \rangle\rangle =
\langle\langle  \lambda^{-1}\H \lambda, \M  \rangle\rangle
\subset \langle\langle \H+ \H^*, \M\rangle\rangle = \{0\}$ by
Lemma 3.1.
{\em Q.E.D.}

We are now ready to prove our main auxiliary statement.

\smallskip
\noindent
{\bf Theorem 3.4.} {\em
Let us adopt the assumptions of Lemma 3.1,
and consider a submanifold $\check H^*\subset K^*$ defined locally
by the constraints $\xi_{M^i}=0$, where the functions $\xi_{M^i}$
are associated by (\ref{3.12}) with a basis $\{M^i\}$ of $\M$.
Then the PBs of the constraints are given by
\be
C^{ij}(\lambda):= \{\xi_{M^i}, \xi_{M^j}\}_{K^*}(\lambda)=
\langle\langle  (\lambda M^j \lambda^{-1})_{\M^*},
(\lambda M^i \lambda^{-1})_\M \rangle \rangle
\qquad (\lambda\in \check H^*).
\label{3.16}\ee
If the matrix $C^{ij}(\lambda)$ is non-degenerate for
$\lambda \in \check H^*$, then  the Dirac reduction of
$(K^*, \{\ ,\ \}_{K^*})$ yields $(\check H^*, \{\ ,\ \}_{\check H^*})$.
}

\medskip
\noindent
{\em Proof.}
Let the functions $F_n$ and $f_n$ be related by (\ref{3.10}) for $n=1,2$.
The statement of the theorem follows by combining the preceding lemmas
with the standard formula \cite{Dirac} of the Dirac bracket, $\{\ ,\ \}^*$:
$$
\{ F_1, F_2\}^*(\lambda) =
\{ f_1, f_2\}_{K^*}(\lambda) - \sum_{i,j}
\{ f_1, \xi_{M^i}\}_{K^*}(\lambda) (C^{-1}(\lambda))_{ij} \{\xi_{M^j},
f_2\}_{K^*}(\lambda).
$$
The second term vanishes by (\ref{3.14}),
and the first term yields $\{ F_1, F_2\}_{H^*}$ on account of
Lemmas 3.1 and 3.2. {\em Q.E.D.}

\medskip
\noindent
{\em Remark 3.5.}
It is clear from the proof of Theorem 3.4 that the above used
assumptions are not only sufficient, but also necessary for
the desired Dirac reduction to work.
For example, the assumption (\ref{3.5}) is crucial in the
proof of Lemma 3.1 on which Theorem 3.4 relies;
we were led to this assumption in the $R=0$ case studied in \cite{FGP}, too.
The same assumptions appear in the construction of
PL dynamical $r$-matrices given in \cite{EEM}.
In a sense, Dirac reductions provides (for us) an
explanation of these assumptions.

\subsection{PL dynamical $r$-matrices from Dirac reduction}

By using the framework developed so far, the
following result is essentially obvious.

\smallskip
\noindent
{\bf Theorem 3.6.}
{\em
Consider $(Q(\check K^*), \{\ ,\ \}_{Q(\check K^*)})$ with the PB defined by a
(possibly zero) triangular PL dynamical $r$-matrix $r:\check K^* \rightarrow \G\wedge \G$.
Adopt the assumptions of Lemma 3.1 and suppose that
$C^{ij}(\lambda)$ (\ref{3.16}) gives a non-degenerate matrix function
on a non-empty submanifold $\check H^*\subset H^*$, which is contained
in $\check K^*$.
Then the submanifold $Q(\check H^*) \subset Q(\check K^*)$
is defined by second class constraints, and the resulting Dirac bracket
is of the type $\{\ ,\ \}_{Q(\check H^*)}$ with the reduced $r$-matrix
\be
r^*(\lambda) = r(\lambda) + \rho(\lambda)  \qquad
 (\lambda \in \check H^*),
\label{3.18}\ee
where
$\rho: \check H^* \rightarrow \M\wedge \M\subset
\K\wedge\K\subset \G\wedge \G$ is given by
\be
\rho(\lambda)=
\sum_{i,j} (C^{-1}(\lambda))_{ij} (\lambda M^i \lambda^{-1})_\M
\otimes (\lambda M^j \lambda^{-1})_\M,
\qquad  \forall \lambda\in \check H^*.
\label{3.19}\ee
Here, $\{ M^i\}$ is a basis of $\M$ (\ref{3.3}) and
$\Ad_\lambda M^i=\lambda M^i \lambda^{-1}$ is defined using
the double $D(K,K^*)$.}

\noindent
{\em Proof.}  One can easily calculate the Dirac bracket similarly
to the proof of Theorem 3.4.
{\em Q.E.D.}

\medskip
\noindent
{\bf Corollary 3.7.} {(by Proposition 2.1).}
{\em
Let $r: \check K^* \rightarrow \G\wedge \G$ be a (possibly zero)
triangular PL dynamical
$r$-matrix for $K\subseteq G$.
Suppose that the assumptions of Lemma 3.1 hold and $C^{ij}$ (\ref{3.16})
defines a non-degenerate matrix function on a non-empty submanifold, $\check H^*$,  of $H^*\cap \check K^*$.
Then $r^*: \check H^* \rightarrow \G\wedge \G$ (\ref{3.18}) gives a triangular PL dynamical
$r$-matrix for $H\subset G$.}

\medskip

If the bases $\{M_i\}\subset \M^*$ and $\{M^i\}\subset \M$ are in duality,
then so are the bases
$\{\lambda M_i \lambda^{-1}\}\subset \M^*$ and
$\{(\lambda M^i\lambda^{-1})_\M\}\subset \M$ for
any  $\lambda\in \check H^*$.
By the invertibility of $C^{ij}(\lambda)$,
$\{(\lambda M^i\lambda^{-1})_{\M^*}\}\subset \M^*$ forms a basis, too.
It follows that for any base element $M_i\in \M^*$ and  $\lambda\in \check H^*$
there exists a unique element $N_i(\lambda)$ that satisfies
\be
\lambda^{-1} M_i \lambda = (\lambda^{-1} N_i(\lambda) \lambda)_{\M^*},
\qquad
N_i(\lambda) \in \M.
\label{3.20}\ee

\smallskip
\noindent
{\bf Lemma 3.8.}
{\em By using $N_i(\lambda)$ (\ref{3.20}),
the triangular PL $r$-matrix in (\ref{3.19}) can be written as
\be
\rho(\lambda) = -\sum_i N_i(\lambda) \otimes M^i=
\sum_i M^i\otimes N_i(\lambda)
\qquad
\forall \lambda \in \check H^*.
\label{3.21}\ee
}

\medskip
\noindent
{\em Proof.}
We have to show that the operator
 $\hat \rho(\lambda) \in \mathrm{End}(\M^*,\M)$, defined by
$$
\hat \rho(\lambda)(M_k)= \sum_{i,j}
(C^{-1}(\lambda))_{ij} (\lambda M^i \lambda^{-1})_\M
\langle M_k, (\lambda M^j \lambda^{-1})_\M \rangle,
$$
satisfies $\hat \rho(\lambda) (M_k)=-N_k(\lambda)$.
By the definition of $N_k(\lambda)$ and the invariance
of the scalar product of ${\cal D}(\K,\K^*)$, we have
\bea
&&\hat \rho(\lambda)(M_k)
=
\sum_{i,j} (C^{-1}(\lambda))_{ij} (\lambda M^i \lambda^{-1})_\M
\langle \langle (\lambda^{-1} N_k(\lambda)  \lambda)_{\M^*},
 M^j \rangle\rangle \nonumber\\
&&\quad =
\sum_{i,j}
(C^{-1}(\lambda))_{ij} (\lambda M^i \lambda^{-1})_\M
\langle \langle  N_k(\lambda),
\lambda M^j \lambda^{-1})_{\M^*}\rangle\rangle \nonumber\\
&&\quad =
\sum_{i,j,l}
(C^{-1}(\lambda))_{ij} (\lambda M^i \lambda^{-1})_\M
\langle \langle  N_k(\lambda),
\lambda M_l \lambda^{-1} \rangle\rangle C^{lj}(\lambda)
 = - N_k(\lambda),
\nonumber\eea
as required. {\em Q.E.D.}

\smallskip
\noindent
{\em Remark 3.9.} The dynamical $r$-matrix $\rho$ in (\ref{3.19})
is the same as $\sigma_\H^\G$ found in \cite{EEM}.
In order to verify this, note that formula (\ref{3.21}) implies
the identity
$$
\langle \langle (\lambda^{-1} u \lambda)_\M,
\lambda^{-1} v \lambda\rangle\rangle =
\sum_i \langle \langle (\lambda^{-1} u \lambda)_\M,
\lambda^{-1} M^i \lambda  \rangle\rangle \,
\langle \langle (\lambda^{-1} v \lambda)_\M,
\lambda^{-1} N_i(\lambda) \lambda \rangle\rangle
$$
for all $u,v\in \M$, $\lambda\in \check H^*$.
According to \cite{EEM} (property 1
above Theorem 2.2) this identity
characterizes $\sigma_\H^\G$ uniquely, if $\sigma_H^\G$
is written in the form (\ref{3.21}) {\em with some}  $N_i(\lambda)$.
The arguments that led to our Corollary 3.7 appear (for us)
more enlightening than
the direct proof of Theorem 2.2 in \cite{EEM}, which states
that $\sigma_\H^\G$ is a triangular PL dynamical $r$-matrix.

\section{Discussion}
\setcounter{equation}{0}

It is important to note that
the applicability of the Dirac reduction method is not restricted
to the triangular case.
In fact \cite{FM1,EEM}, an arbitrary PL dynamical $r$-matrix
$r: \check K^* \rightarrow \G\wedge \G$ encodes a PB on
the manifold
\be
P=P(\check K^*):= \check K^* \times G \times \check K^* =
\{ (\tilde \kappa, g, \hat \kappa)\}.
\label{4.1}\ee
For admissible
functions $f\in \F(\check K^*)$ and $\phi\in \F(G)$ one introduces
$\hat f, \tilde f\in \F(P)$ and $\phi'\in \F(P)$
by
$\hat f(\tilde \kappa, g, \hat\kappa) = f(\hat \kappa)$,
$\tilde f(\tilde \kappa, g, \hat\kappa) = f(\tilde \kappa)$,
$\phi'(\tilde \kappa, g, \hat\kappa) = \phi(g)$.
One then postulates a bracket on the functions on $P$ by the ansatz
\bea
&\{ \hat f_1, \hat f_2\}_{P}(\tilde\kappa, g, \hat \kappa)=
\{ f_1, f_2 \}_{\check K^*}(\hat \kappa),
\qquad
\{ \tilde f_1, \tilde f_2\}_{P}(\tilde\kappa, g, \hat\kappa)=
-\{ f_1, f_2 \}_{\check K^*}(\tilde \kappa),
\nonumber\\
&\{ \phi', \hat f\}_{P}(\tilde\kappa, g, \hat\kappa)=
\langle \nabla'_g \phi, \nabla_{\hat\kappa} f\rangle,
\qquad
\{ \phi', \tilde f\}_{P}(\tilde\kappa, g, \hat\kappa)=
\langle \nabla_g \phi, \nabla_{\tilde \kappa} f\rangle,
\label{4.2}\\
&\{ \phi'_1, \phi'_2 \}_{P}(\tilde \kappa, g,\hat\kappa)=
\langle \nabla'_g \phi_1 \otimes \nabla'_g \phi_2, R+r(\hat \kappa) \rangle
- \langle \nabla_g \phi_1 \otimes \nabla_g \phi_2, R+
r(\tilde \kappa) \rangle,
\nonumber\eea
together with $\{ \hat f_1, \tilde f_2\}_{P}=0$, where
$f, f_i \in \F(\check K^*)$, $\phi, \phi_i\in \F(G)$,
$\langle\ ,\ \rangle$ denotes the canonical pairing
between elements of $\G^*$ and $\G$, and $R$ is the chosen
constant $r$-matrix (\ref{1.2}).
The ansatz (\ref{4.2}) defines a PB if and only if the PL-CDYBE
(\ref{1.3}) and the equivariance condition
(\ref{2.10}) are valid for $r$.
It is clear that under the assumptions of Theorem 3.4
the Dirac reduction of $(P(\check K^*), \{\ ,\ \}_{P(\check K^*)})$
yields $(P(\check H^*), \{\ ,\ \}_{P(\check H^*)})$, and the
accompanying reduction of the PL dynamical $r$-matrix is given
by the same formula (\ref{3.18}), (\ref{3.19})  as in the triangular case.
The content of this  statement is precisely
the `composition theorem' (Theorem 2.7) of \cite{EEM}.
(Note also that ${\cal I}_{R,r}={\cal I}_{R,r^*}$ is easily checked by
using (\ref{3.21}).)

If $R=0$, then the construction of dynamical $r$-matrices by
Dirac reduction described above
specializes to the construction given in \cite{FGP}.
This provides us with examples in the  case of an Abelian $G^*$.
For non-Abelian $G^*$ we do not know examples
that are essentially different from those mentioned in \cite{EEM}.
If $R$ is the standard (Drinfeld-Jimbo)
{\em factorisable} $r$-matrix on a simple Lie
algebra, then one can apply the reduction by taking
 $\K=\G$ and taking $\H$ to be a Levi (regular reductive)
subalgebra of $\G$.
Thus Corollary 3.7 yields triangular PL dynamical $r$ matrices for the
Levi subgroups of $G$.
The composition theorem can also be applied by taking
the $r_{\mathrm{BFP}}$ solution
\cite{BFP} of the PL-CDYBE for $K=G$ as
 the starting point \cite{EEM}.
Although not mentioned in \cite{EEM}, the same family of
 examples is available
in the compact case as well, where a simple compact Lie group $G$
is equipped with its
standard PL structure and $H\subset K=G$ is a regular reductive
subgroup.
(See also \cite{FM2} for
a description of $r_{\mathrm{BFP}}$ in PL terms.)

Incidentally,
the Dirac reductions of $r_{\mathrm{BFP}}$ just alluded to
can be seen as exchange $r$-matrices in the
Wess-Zumino-Novikov-Witten model,
obtained there by restricting the monodromy matrix to a regular
reductive subgroup of $G$, i.e., by performing the corresponding
Dirac reduction of the chiral WZNW PB defined by
$r_{\mathrm{BFP}}$ \cite{BFP}.
The closely related trigonometric PL $r$-matrices of \cite{M}
can also be associated with suitable  PL symmetries on the chiral
WZNW phase space with restricted  monodromy.

It appears an interesting open question whether one can relate
the PL dynamical $r$-matrices
to finite dimensional integrable systems by suitable
extension of the constructions in \cite{LiXu,Li} to the PL case.
Our basic idea for such generalization is to apply Hamiltonian reduction to
$(P(\check K^*), \{\ ,\ \}_{P(\check K^*)})$ by using the PL action
of $K$ generated by the PL momentum map $\Lambda: P\to K^*$ given by
\be
\Lambda: (\tilde \kappa, g, \hat \kappa)\mapsto
\tilde \kappa \hat \kappa^{-1}.
\label{4.3}\ee
In analogy with the constructions in \cite{LiXu,Li},
a relevant reduction of $P$ should be defined by setting
the momentum map $\Lambda$ to unity;
in other words by imposing the first class constraints
$\hat \kappa = \tilde \kappa$.
However, we have not yet investigated how to obtain commuting
Hamiltonians on the reduced phase space in this context.
Naively, one expects to obtain such Hamiltonians from
 the functions
of the form $\phi'$ where $\phi$ is a central function on $G$,
but further work is required to see if this idea can really work or not.

It could be also
interesting to develop the quantum version of our Dirac reduction
algorithm.
This may simplify the quantization
of  various dynamical $r$-matrices \cite{EE,EEM}.

\bigskip
\bigskip
{\small
\noindent{\bf Acknowledgements.}
This work was supported in part by the Hungarian
Scientific Research Fund (OTKA) under grants
T034170, T043159 and M036803.
Thanks are due to J. Balog and
I. Marshall for useful comments on the remark
communicated here.
}


\newpage
\begin{thebibliography}{EnGH}

\bibitem{BFP}
J. Balog, L.  Feh\'er and  L. Palla,
{\sl Chiral extensions of the WZNW phase space,
Poisson-Lie symmetries and groupoids},
Nucl. Phys. B {\bf 568} (2000), 503-542 (hep-th/9910046).

\bibitem{FM1}
L. Feh\'er and I. Marshall,
{\sl On a Poisson-Lie analogue of the classical
dynamical Yang-Baxter equation for self-dual Lie algebras},
Lett. Math. Phys. {\bf 62} (2002), 51-62
(math.QA/0208159).

\bibitem{FM2}
L. Feh\'er and I. Marshall, {\sl The non-Abelian momentum map for
Poisson-Lie symmetries on the chiral WZNW phase space},
math.QA/0401226.

\bibitem{EEM}
B. Enriquez, P. Etingof and I. Marshall,
{\sl Quantization of some Poisson-Lie dynamical $r$-matrices and
 Poisson homogeneous spaces}, math.QA/0403283.

\bibitem{DM}
J. Donin and A. Mudrov,
{\sl Dynamical Yang-Baxter equation and quantum vector bundles},
math.QA/0306028.

\bibitem{M}
A. Mudrov,
{\sl Trigonometric dynamical $r$-matrices over Poisson Lie base},
math.QA/0403207.

\bibitem{EV}  P. Etingof and A. Varchenko,
{\sl Geometry and classification of solutions of the classical
dynamical Yang-Baxter equation},
Commun. Math. Phys. {\bf 192} (1998), 77-129 (q-alg/9703040).

\bibitem{FGP}
L. Feh\'er,  A. G\'abor and  B.G. Pusztai,
{\sl On dynamical $r$-matrices obtained from Dirac reduction
and their generalizations to affine Lie algebras},
J. Phys. A {\bf 34} (2001), 7235-7248 (math-ph/0105047).

\bibitem{EE}
B. Enriquez and P. Etingof,
{\sl Quantization of classical dynamical $r$-matrices with
nonabelian base}, math.QA/0311224.

\bibitem{LW}
J.-H. Lu and A. Weinstein,
{\sl Poisson Lie groups, dressing transformations and Bruhat decompositions},
J. Diff. Geom. {\bf 31} (1990), 501-526.

\bibitem{Dirac}
P.A.M. Dirac,
Lectures on Quantum Mechanics, Yeshiva University Press, 1964.

\bibitem{STS}
M.A. Semenov-Tian-Shansky,
{\sl Dressing transformations and Poisson group actions},
Publ. RIMS {\bf 21} (1985), 1237-1260.

\bibitem{LiXu}
L.C. Li and P. Xu,
{\sl Integrable spin Calogero-Moser systems},
Commun. Math. Phys. {\bf 231} (2002), 257-286 (math.QA/0105162).

\bibitem{Li}
L.C. Li,
{\sl Coboundary dynamical Poisson-Li groupoids and integrable
systems}, Int. Math. Res. Not. {\bf 2003 No. 51} (2003),  2725-2746.

\end{thebibliography}
\end{document}